\documentclass{colt2015} 

\usepackage{algorithm,algorithmic}
\usepackage{amssymb, multirow, paralist, color}
\newtheorem{thm}{Theorem}

\newtheorem{ass}[thm]{Assumption}
\usepackage{enumitem}
\usepackage{textcase}
\usepackage[T1]{fontenc}

\def \R {\mathbb{R}}

\def \v {\mathbf{v}}

\def \x {\mathbf{x}}

\def \x {\mathbf{x}}

\def \e {\mathbf{e}}

\def \z {\mathbf{z}}

\def \y {\mathbf{y}}
\def \u {\mathbf{u}}
\def \uh {\widehat{\u}}

\def \xh {\widehat{\x}}

\def \y {\mathbf{y}}

\def \x {\mathbf{x}}

\def \z {\mathbf{z}}
\def \u {\mathbf{u}}

\def \R {\mathbb{R}}

\def \neon {\textsc{Neon}}
\def \neonp {\textsc{Neon}$^+$}

\def \v {\mathbf{v}}

\def \xh {\widehat{\x}}

\begin{document}

\title[AG methods for Extracting Negative Curvature From Noise]{\textsc{Neon$^+$}: Accelerated Gradient Methods for Extracting Negative Curvature for Non-Convex Optimization\thanks{The main result of extracting negative curvature by \textsc{Neon$^+$} in this manuscript  is merged into our earlier manuscript ``First-order Stochastic Algorithms for Escaping From Saddle Points in Almost Linear Time''~\citep{NEON17}. }}
 \author{\Name{Yi Xu}$^\natural$\Email{yi-xu@uiowa.edu}\\
 \Name{Rong Jin}$^\dagger$ \Email{jinrong.jr@ailibaba-inc.com}\\
\Name{Tianbao Yang}$^\natural$\Email{tianbao-yang@uiowa.edu}\\
\addr$^\natural$ Department of Computer Science, The University of Iowa, Iowa City, IA 52242 \\
\addr$^\dagger$Alibaba Group, Bellevue, WA 98004 \\
}

\maketitle
\vspace*{-0.5in}
\begin{center}{First version: December 02, 2017}\end{center}
\vspace*{0.2in}

\begin{abstract}
Accelerated gradient (AG) methods are breakthroughs in convex optimization, improving the convergence rate of the gradient descent method for optimization with smooth functions. However, the analysis of AG methods  for non-convex optimization is still limited.  It remains an open question whether AG methods from convex optimization can accelerate the convergence of  the gradient descent method for finding local minimum of non-convex optimization problems. This paper provides an affirmative answer to this question.  In particular, we analyze Nesterov's Accelerated Gradient method for extracting the negative curvature from random noise, which is central to escaping from saddle points. By leveraging  the proposed NAG methods for extracting the negative curvature, we present a new AG algorithm with double loops for non-convex optimization
, which converges to second-order stationary point $\x$ such that $\|\nabla f(\x)\|\leq \epsilon$ and $\nabla^2 f(\x)\geq -\sqrt{\epsilon} I$ with $\widetilde O(1/\epsilon^{1.75})$ iteration complexity, improving that of gradient descent method  by a factor of $\epsilon^{-0.25}$ and matching the best iteration complexity of second-order  Hessian-free methods for non-convex optimization. 
\end{abstract}

\section{Introduction}
We consider the following optimization problem in this paper: 
\begin{align}\label{eqn:opt}
\min_{\x\in\R^d}f(\x),
\end{align}
where $f(\x)$ is a twice differentiable non-convex smooth function, whose Hessian is Lipschitz continuous. Recently, this problem has received increasing interests in the machine learning community due to that many learning problems are non-convex~\cite{}. A renowned  method in the machine learning community is gradient descent (GD) method, which updates the solution according to the following equation: 
\begin{align*}
\x_{t+1} = \x_t - \eta\nabla f(\x_t),
\end{align*}
where $\eta$ is a constant step size. It is not difficult to show that GD converges to an $\epsilon$-critical point $\x$, i.e., $\|\nabla f(\x)\|\leq \epsilon$, with an iteration complexity of $O(1/\epsilon^2)$.
However, such a critical point could be a saddle point, which could be far from a local minimum and can harm the performance of prediction in machine learning. 

To address this issue, one solution is to design an optimization algorithm that can converge to an $(\epsilon, \gamma)$-second-order stationary point (SSP) $\x$ such that 
\begin{align*}
\|\nabla f(\x)\|\leq \epsilon, \quad \nabla^2 f(\x)\geq -\gamma I. 
\end{align*}
When the objective function is non-degenerate (i.e., the Hessian matrix at all saddle points has negative eigen-values),  an $(\epsilon, \gamma)$-SSP $\x$ is guaranteed to be close to a local minimum and even a global minimum for certain problems~\cite{}.

Despite the popularity of GD, second-order methods have emerged to provide second-order convergence guarantee for non-convex optimization, which utilize the Hessian matrix or the Hessian-vector product for updating the solution. Starting from the seminal work by~\citep{nesterov2006cubic}, a wave of studies have been devoted to designing efficient second-order  optimization algorithms with fast convergence to a SSP~\citep{DBLP:conf/stoc/AgarwalZBHM17,DBLP:journals/corr/CarmonDHS16,peng16inexacthessian,Cartis2011,Cartis2011b,clement17,DBLP:journals/corr/noisynegative,DBLP:journals/corr/CarmonD16}. A state-of-the-art result for second-order optimization algorithms is to use Hessian-vector products  for finding an $(\epsilon, \sqrt{\epsilon})$-SSP with  an $\widetilde O(T_h/\epsilon^{1.75})$ time complexity where $T_h$ denotes the runtime of the Hessian-vector product.

In spite of the theoretical promise of second-order optimization algorithms,  first-order (or gradient-based) algorithms are still the first choice in practice due to their simplicity.   A recent breakthrough for gradient-based non-convex optimization methods is due to~\cite{DBLP:conf/icml/Jin0NKJ17}, who proposed a  gradient-based method converging to an $(\epsilon, \sqrt{\epsilon})$-SSP with an iteration complexity of $\widetilde O(1/\epsilon^2)$. Although it promotes  GD for finding a SSP, it  is  still worse by a factor  of $\epsilon^{-0.25}$ than the state-of-the-art second-order optimization algorithms.  Given the dramatic  success of accelerated gradient methods for convex optimization, an interesting question is:
\\
\textit{Can we use accelerated gradient methods from convex optimization to accelerate the convergence  of non-convex optimization for finding a SSP?  }
\\
This paper gives an affirmative answer to this question. Our main contribution is summarized below: 
\begin{itemize}
\item We analyze 
Nesterov's Accelerated Gradient (NAG) method~\citep{opac-b1104789} for extracting the negative curvature of a Hessian matrix, which is central to algorithms escaping from (non-degenerate) saddle points for non-convex optimization. We refer to the proposed procedure for extracting the negative curvature of a Hessian matrix as \neonp since it is an accelerated variant of its predecessor~\neon~\citep{NEON17}. 
\item By combining the proposed AG methods for extracting the negative curvature of a Hessian matrix and an existing AG method for minimizing a regularized almost-convex function~ \citep{DBLP:journals/corr/CarmonDHS16}, we present a new AG algorithm with double  loops (dubbed \textbf{NEAG}) for finding a SSP to a general non-convex optimization problem in~(\ref{eqn:opt}). 
\item The proposed NEAG algorithm enjoys an iteration complexity of $\widetilde O(1/\epsilon^{1.75})$ for finding an $(\epsilon, \sqrt{\epsilon})$-SSP, matching the state-of-the-art result of  second-order Hessian-free methods. 
\end{itemize}

\section{Related Work}
Although there are extensive studies about AG methods for convex optimization, the analysis of AG for non-convex optimization is still limited.  \cite{DBLP:journals/mp/GhadimiL16} analyzed a variant of AG for minimizing non-convex smooth functions. However, its rate of convergence to a critical point is the same as standard GD method.  \cite{Li:2015:APG:2969239.2969282} analyzed variants of accelerated proximal gradient (APG) methods for minimizing  a family of non-convex functions consisting of a smooth component and a non-smooth component. For general functions in this family, they only proved the asymptotic convergence to a critical point. Asymptotic results with  explicit  convergence rates are established for functions that satisfy the Kurdyka -\ Lojasiewicz (KL) property. \cite{yangnonconvexmo} analyzed two variants of AG methods in a stochastic setting for non-convex optimization under a unified framework. Again, their convergence analysis are only for finding critical points and the convergence rates of the analyzed two AG methods in the stochastic setting is the same as stochastic gradient method. 

Recently, \cite{corrACGWright} proved that the Polyak's heavy-ball method method does not converge to critical points that do not satisfy second-order necessary
conditions. They also analyzed the divergence rate of two accelerated gradient methods (including Polyak's heavy-ball method and Nesterov's AG method) from a (non-degenerate) saddle point of a non-convex {quadratic function}, showing  that both methods can diverge from this point more rapidly than GD. Nevertheless, their analysis does not provide any guarantee on permanent escaping from a (non-degenerate) saddle point. In addition,  no explicit convergence rate of AG methods to a SSP was established for a general non-convex optimization problem~(\ref{eqn:opt}). \cite{DBLP:journals/corr/CarmonDHS16} analyzed an AG method for minimizing an almost-convex function, which achieves faster convergence to a critical point than GD method. They also developed an accelerated method for finding a SSP of a general non-convex optimization problem by combining the AG method for minimizing an almost-convex function and a second-order method (e.g., the Lanczos method) for extracting the negative curvature. For finding an $(\epsilon, \sqrt{\epsilon})$-SSP, their accelerated method achieves the best iteration complexity of $\widetilde O(1/\epsilon^{1.75})$ among existing second-order Hessian-free methods. However, they do not address the question raised before, i.e., finding a SSP with an AG method. \cite{DBLP:conf/icml/CarmonDHS17} proposed the first AG algorithm with provable acceleration over GD for non-convex optimization. Their algorithm can converge to a first-order stationary point with $\widetilde O(1/\epsilon^{1.75})$ iterations. However, it is not clear whether their method can guarantee finding a SSP with an iteration complexity of $\widetilde O(1/\epsilon^{1.75})$.

A fundamental concern in the design of non-convex optimization algorithms for finding a SSP is how to escape from saddle points.  The proposed first-order method dubbed \textsc{Neon$^+$} addresses this concern by extracting the negative curvature from a Hessian matrix with negative eigen-values.  It is inspired by a recent work~\citep{NEON17}, which is the first work that develops a first-order method (named \textsc{Neon}) for extracting negative curvature from a Hessian matrix with negative eigen-values. Their method is a gradient descent method and  suffers from an iteration complexity of $\widetilde O(1/\gamma)$ for finding a negative curvature for a Hessian matrix whose minimum eigen-value is less than $-\gamma<0$. In contrast, the proposed \textsc{Neon}$^+$ is based on PHB or NAG and  improves the iteration complexity of \textsc{Neon} to $\widetilde O(1/\sqrt{\gamma})$. By utilizing \textsc{Neon}$^+$ in the framework developed by~\cite{DBLP:journals/corr/CarmonDHS16}, we obtain an AG algorithm for finding an $(\epsilon, \sqrt{\epsilon})$-SSP with an iteration complexity of $\widetilde O(1/\epsilon^{1.75})$. As a byproduct, \neonp~can be also leveraged in stochastic non-convex optimization to accelerate the convergence for finding a SSP. 

It is worth mentioning that, a recent work~\citep{NEON2} also develops an improved variant of \textsc{Neon} named \textsc{Neon2}$^\text{det}$ for finding a negative curvature with an iteration complexity of $\widetilde O(1/\sqrt{\gamma})$. We emphasize three differences between  \textsc{Neon}$^+$ and  \textsc{Neon2}$^\text{det}$: (i) \textsc{Neon}$^+$ is based on PHB and NAG, which are more familiar to the machine learning and optimization community; while  \textsc{Neon2}$^\text{det}$ is based on Chebyshev approximation theory; (ii) our analysis of \textsc{Neon}$^+$  is elementary and self-contained; while the analysis of \textsc{Neon2}$^\text{det}$ relies on the stability analysis of Chebyshev  polynomials; (iii)  \textsc{Neon2}$^\text{det}$  terminates when the resulting vector is within an Euclidean ball with a small radius inversely proportional to a power of the dimensionality, which might cause numerical issues in practice for high-dimensional problems; in contrast  \textsc{Neon}$^+$  terminates when the resulting vector is within an Euclidean ball with a  radius almost independent of the dimensionality, rendering it much more viable for high dimensional problems.

It was also brought to our attention that  when we prepare this manuscript,  an independent work by~\citep{AGNON} also analyzed  Nesterov's AG method for non-convex optimization. For finding a $(\epsilon, \sqrt{\epsilon}$)-SSP, the algorithm~\citep{AGNON} and in the present work enjoy the same iteration complexity of $\widetilde O(1/\epsilon^{1.75})$. The differences between the two work are (i) we focus our analysis on 
NAG for extracting negative curvature, in contrast they directly analyze NAG for non-convex optimization; (ii) our AG algorithm has a nested loop, while their algorithm is a single loop. 


\section{Preliminaries}
In this section, we present some preliminaries, including some notations, accelerated gradient methods for convex optimization, and the idea of  \textsc{Neon} for non-convex optimization.

Denote by $\|\cdot\|$ the Euclidean norm of a vector and $\|\cdot\|_2$ the spectral n orm of  a matrix. Let $\lambda_{\min}(\cdot)$ denote the minimum eigen-value of a matrix. We make standard assumptions regarding~(\ref{eqn:opt}) in order to find a SSP. 
\begin{ass}
\begin{enumerate}
\item $f(\x)$ has $L_1$-Lipschitz continuous gradient and $L_2$-Lipschitz continuous Hessian, i.e., 
\begin{align}
\|\nabla f(\x) - \nabla f(\y)\|\leq L_1\|\x - \y\|, \quad \|\nabla^2 f(\x) - \nabla^2 f(\y)\|_2\leq L_2\|\x - \y\|
\end{align}
\item given an initial point $\x_0$, assume that there exists $0<\Delta<\infty$ such that $f(\x_0) - \min_{\x\in\R^d}f(\x)\leq \Delta$.
\end{enumerate}
\end{ass}
A function $g(\x)$ is $\sigma$-strongly convex ($\sigma>0$) if for all $\x,\y$ it holds that
\begin{align*}
g(\x)\geq g(\y) + \nabla g(\y)^{\top}(\x - \y)  + \frac{\sigma}{2}\|\x - \y\|^2, \forall \x, \y
\end{align*}
If the above inequality holds for $\sigma<0$, then $g(\x)$ is called $(-\sigma)$-almost convex.

\subsection{Accelerated Gradient Methods from Convex Optimization}
For minimizing a smooth and (strongly) convex function $g(\x)$, AG methods have been proposed with faster convergence rate than GD. Next, we present a variant of AG  for minimizing a $L_1$-smooth and $\sigma_1$-strongly convex function  since it will be used in our development. A famous  variant of AG is Nesterov's AG (NAG) method~\citep{opac-b1104789}, whose update is given by Step 6 \&7 in Algorithm~\ref{alg:agd}. The iteration complexity of Algorithm~\ref{alg:agd} for minimizing a $L_1$-smooth and $\sigma_1$-strongly convex function is given by $O(\sqrt{L_1/\sigma_1}\log(1/\epsilon))$. 


\begin{algorithm}[t]
\caption{An AG method for minimizing a Smooth and Strongly Convex function: \newline
AG-SSC$(g,\y_1,\epsilon,L_1,\sigma_1)$
}\label{alg:agd}
\begin{algorithmic}[1]
\STATE Set $\kappa=L_1/\sigma_1$, $\z_1=\y_1$, $\zeta =  \frac{\sqrt{\kappa}-1}{\sqrt{\kappa}+1}$
\FOR{$j=1,2,\ldots$}
\IF{$\|\nabla g(\y_j)\|\leq\epsilon$}
\RETURN $\y_j$
\ENDIF
\STATE $\y_{j+1}=\z_j-\frac{1}{L_1}\nabla g(\z_j)$,  
\STATE $\z_{j+1}=\y_{j+1} +\zeta (\y_{j+1}-\y_j)$
\ENDFOR
\end{algorithmic}
\end{algorithm}

\subsection{\neon~for Non-Convex Optimization}
\neon~\citep{NEON17}  is a first-order procedure for extracting the negative curvature from a Hessian matrix $\nabla^2 f(\x)$ at a point $\x$. In particular, it starts from a random noise vector $\u_0$ and iteratively updates $\u_\tau$ according to the following equation: 
\begin{align}\label{eqn:neon}
\u_{\tau+1} = \u_\tau - \eta (\nabla f(\x+\u_\tau) - \nabla f(\x)), 
\end{align}
and it terminates when $f(\x+\u) - f(\x) - \nabla f(\x)^{\top}\u$ is sufficiently small given that $\u_\tau$ resides in an Euclidean ball with a proper radius, or it cannot make $f(\x+\u) - f(\x) - \nabla f(\x)^{\top}\u$ sufficiently small after a certain number of iterations. In~\citep{NEON17}, \neon~is motivated by using a noisy Power method to compute the negative curvature. In particular, the Lipschitz Hessian condition implies that 
\begin{align*}
\|\nabla f(\x+\u_\tau) - \nabla f(\x) - \nabla^2 f(\x)\u_\tau\|\leq \frac{L_2}{2}\|\u\|^2.
\end{align*}
Therefore, the sequence in~(\ref{eqn:neon}) approximates another sequence $\u'_{\tau+1}= (I - \eta \nabla^2 f(\x))\u_\tau'$ when $\|\u_\tau\|$ is sufficiently small, which is exactly the sequence generated by applying the Power method to compute the leading eigen-pair of $I - \eta \nabla^2 f(\x)$ corresponding to the minimum eigen-pair of $\nabla^2 f(\x)$. 

On the other hand, we can also consider the sequence~(\ref{eqn:neon}) as an application of GD to the following objective function:
\begin{align}\label{eqn:objn}
\hat f_\x(\u) = f(\x+\u) - f(\x) - \nabla f(\x)^{\top}\u. 
\end{align}
Sometimes write $\hat f_\x(\u) = \hat f(\u)$, where the dependent $\x$ should be clear from the context. 
By the Lipschitz continuous Hessian condition, we have that
\begin{align*}
\frac{1}{2}\u^{\top}\nabla^2 f(\x)\u - \frac{L_3}{6}\|\u\|^3 \leq \hat f(\u) \leq \frac{1}{2}\u^{\top}\nabla^2 f(\x)\u + \frac{L_3}{6}\|\u\|^3.
\end{align*}
It implies that if $\hat f(\u)$ is sufficiently less than zero and $\|\u\|$ is not too large, then $\frac{\u^{\top}\nabla^2 f(\x)\u}{\|\u\|^2}$ will be sufficiently less than zero.

A key result in~\citep{NEON17} is that  when $\lambda_{\min}(\nabla^2 f(\x))\leq - \gamma $ \neon~can find a negative curvature direction $\u$ such that $\frac{\u^{\top}\nabla^2 f(\x)\u}{\|\u\|^2}\leq -\widetilde\Omega(\gamma)$ with an $\widetilde O(1/\gamma)$ number of iterations. 
The main contribution of this paper is to show that both PHB and NAG with an appropriate momentum constant $\zeta$ when applied the function $\hat f(\u)$ can find a negative curvature much faster than GD.

\section{\textsc{Neon}$^+$: Accelerated Gradient Methods for Extracting Negative Curvature}

In this section, we will analyze 
NAG methods. We first present the updates of \neonp and discuss the underlying intuition why 
NAG can find a negative curvature faster than GD.

The updates of NAG method applied to the function $\hat f(\u)$ at a given point $\x$ is given by 
\begin{equation}\label{eqn:neonpn}
\begin{aligned}
\y_{\tau+1}& = \u_{\tau} - \eta \nabla \hat f(\u_{\tau})  \\
\u_{\tau+1}& =\y_{\tau+1} +  \zeta(\y_{\tau+1} - \y_{\tau})
\end{aligned}
\end{equation}
where the term $\zeta(\y_{\tau+1} - \y_{\tau})$ is the momentum term, and $\zeta\in(0, 1)$ is the momentum parameter.  
The proposed algorithm based on the NAG method (referred to as NEON$^+$) for extracting NC of a Hessian matrix $\nabla^2 f(\x)$ is presented in Algorithm~\ref{alg:neonp},  where  $$\Delta_\x(\y_\tau,\u_\tau)=\hat f_\x(\y_\tau)- \hat f_\x(\u_\tau) -  \nabla \hat f_\x(\u_\tau)^{\top}(\y_\tau - \u_\tau),$$ and NCFind is a procedure that returns a NC by searching over the history $\y_{0:\tau}, \u_{0:\tau}$ shown in Algorithm~\ref{alg:ncfind}. The condition check in Step 4 is to detect easy cases such that NCFind can easily find a NC in historical solutions without continuing the update. It is notable that NCFind is similar to a procedure called Negative Curvature Exploitation (NCE) in~\citep{AGNON}. However, the difference is that NCFind is tailored to finding a negative curvature, while NCE in~\citep{AGNON} is for ensuring a decrease on a modified objective. 

\begin{algorithm}[t]
\caption{Accelerated Gradient methods for Extracting NC from Noise: \neonp$(f, \x, t, \mathcal F, U,  \zeta, r)$}\label{alg:neonp}
\begin{algorithmic}[1]
\STATE \textbf{Input}:  $f, \x, t, \mathcal F,  U, \zeta, r$
\STATE Generate $\y_0=\u_0$ randomly from the sphere of an Euclidean ball of radius $r$
\FOR{$\tau=0,\ldots, t$}
\IF{$\Delta_\x(\y_\tau, \u_\tau)< - \frac{\gamma}{2}\|\y_\tau - \u_\tau\|^2$}
\STATE {\bf return} $\v=$NCFind($\y_{0:\tau}, \u_{0:\tau}$)
\ENDIF
\STATE compute $(\y_{\tau+1}, \u_{\tau+1})$  by~(\ref{eqn:neonpn})
\ENDFOR
\IF{$\min_{\|\y_{\tau}\|\leq U }\hat f_\x(\y_\tau)  \leq -2\mathcal F$ }
\STATE  let ${\tau'}= \arg\min_{\tau, \|\y_\tau\|\leq  U}\hat f_\x(\y_\tau)$
\STATE  {\bf return} $\y_{\tau'}$
\ELSE 
\STATE {\bf return} $0$ 
\ENDIF
\end{algorithmic}
\end{algorithm}
\begin{algorithm}[t]
\caption{NCFind $(\y_{0:\tau}, \u_{0:\tau})$}\label{alg:ncfind}
\begin{algorithmic}[1]
\IF{$\min_{j=0,\ldots, \tau}\|\y_{j} - \u_j\|\geq \zeta\sqrt{6\eta\mathcal{F}}$}
\STATE {\bf return} $\y_j$, \\where $j = \min\{j': \|\y_{j'}-\u_{j'}\|\geq \zeta\sqrt{6\eta\mathcal{F}}\}$
\ELSE 
\STATE {\bf return} $\y_\tau- \u_{\tau}$
\ENDIF
\end{algorithmic}
\end{algorithm}
%
%
Before presenting our main result and formal analysis, we first present an informal analysis about why \neonp~is faster than \neon.  This analysis is from a perspective of Power method for computing dominating eigen-vectors of a matrix in light of our goal is to compute a negative curvature of a Hessian matrix corresponding to negative eigen-values. In particular, by ignoring the error of $\nabla\hat f(\u_\tau) = \nabla f(\x+\u_\tau) - \nabla f(\x)$ for approximating $H=\nabla^2 f(\x)$, 
the update in~(\ref{eqn:neonpn}) can be written as 
\begin{align}\label{eqn:neonp3}
\uh_{\tau+1} = \left[\begin{array}{cc}(1 + \zeta)(I - \eta H)& -\zeta (I - \eta H)\\ I & 0 \end{array}\right] \uh_{\tau-1}. 
\end{align}
where
\begin{align*}
\uh_{\tau+1} = \left[\begin{array}{cc}\u_{\tau+1}\\ \u_\tau\end{array}\right],
\end{align*}
The above sequence can be considered as an application of the Power method to an augmented matrix:
\begin{align*}
A =  \left[\begin{array}{cc}(1 + \zeta)(I - \eta H)& -\zeta (I - \eta H)\\ I & 0 \end{array}\right].
\end{align*}
According to existing analysis of the Power method for computing top eigen-vectors in the top-$k$ eigen-space of the matrix $A$~\cite{}, the iteration complexity depends on the eigen-gap $\Delta_k$ of $A$  in an order of $\widetilde O(1/\Delta_k)$, where the eigen-gap is defined as the difference between the $k$-th largest  eigen-value of $A$ and the $(k+1)$-th largest eigen-value of $A$. The following result exhibits that by appropriately choosing the momentum constant $\zeta$, the eigen-gap $\Delta_k$ of $A$ matrix involved in NAG method scales as $\sqrt{\gamma}$ for a Hessian matrix $H$ whose negative eigen-values are less than $-\gamma$. 
Following by the similar analysis in~\citep{corrACGWright}, we can state the result formally in the following lemma.
\begin{lemma}\label{lem:1}
Assume the eigen-values of $H$ satisfy $\lambda_1\leq\lambda_2\ldots\leq\lambda_k\leq -\gamma<0\leq \lambda_{k+1}\leq\lambda_d$ and $\eta\leq 1/L_1$ is sufficiently small.  Let $\e_1, \ldots, \e_d$ denote the corresponding eigen-values of $H$. Then the top-$k$ eigen-pairs of $A$ are $(\lambda^h_i(A), \widehat\e_i), i=1,\ldots, k$, where $\widehat\e_i=\left[\begin{array}{cc}\e_i\\ (1/\lambda^h_i(A))\e_i\end{array}\right]$ and
\begin{align*}
\lambda_i^{h}(A) =& \frac{1}{2} [(1+\zeta)(1-\eta\lambda_i) +\sqrt{(1+\zeta)^2(1-\eta\lambda_i)^2-4\zeta(1-\eta\lambda_i)} ].
\end{align*}
By choosing $\zeta = 1 -\sqrt{\eta\gamma}\in(0,1)$, the eigen-gap $\Delta_k$ of $A$ corresponding to its top-$k$  eigen-space is at least $\sqrt{\eta\gamma}/2$.
\end{lemma}
The proof of above lemma can be found in~\citep{NEON17}. We emphasize that the above informal analysis is not enough for proving our main result. 
Nevertheless, a formal analysis yields the following result. 
\begin{thm}\label{thm:main:AGD}
For any $\gamma\in(0,1)$ and a sufficiently small $\delta\in(0,1)$, let $\x$ be a point such that $\lambda_{\min}(\nabla^2 f(\x))\leq -\gamma$. For any constant $\hat c\geq 43$, there  exists a constant $c_{\max}$ that depends on $\hat c$, such that  if \neonp is called with $t = \sqrt{\frac{\hat c\log (dL_1 /(\gamma\delta))}{\eta \gamma}}$,  $\mathcal F=\eta  \gamma^3 L_1L_2^{-2} \log^{-3}(dL_1 /(\gamma\delta))$, $r = \sqrt{\eta }\gamma^2L_1^{-1/2}L_2^{-1}\log^{-2}(dL_1 /(\gamma\delta))$, $U=12\hat c(\sqrt{\eta L_1}\mathcal F/L_2)^{1/3}$, a small constant $\eta\leq c_{\max}/L_1$, and a momentum parameter  $\zeta = 1 - \sqrt{\eta \gamma}$, then  with high probability $1-\delta$ it returns a vector $\u$ such that 
$\frac{\u^{\top}\nabla^2 f(\x)\u}{\|\u\|^2} \leq -\frac{\gamma}{72\hat c^2 \log(dL_1 /(\gamma\delta))} \leq  -\widetilde\Omega(\gamma)$. 
If \neonp returns $\u\neq 0$, then the above inequality must hold; if \neonp returns $0$, we can conclude that $\lambda_{\min}(\nabla^2 f(\x))\geq -\gamma$ with high probability $1- O(\delta)$. 
\end{thm}
{\bf Remark.} The proof of above theorem can be found in~\citep{NEON17}. 

\section{Applications}

\subsection{An AG Algorithm for Non-Convex Optimization}
Built upon \neonp~and an existing framework proposed in~\citep{DBLP:journals/corr/CarmonDHS16}, we can obtain an accelerated algorithm using accelerated gradient methods from convex optimization for non-convex optimization. The proposed AG algorithm is based on two building blocks: (i) an AG method based on \neonp~conducting negative curvature descent for driving the solution to reach a point where the objective function $f(\cdot)$ is locally almost-convex (i.e., the Hessian matrix at the point has all eigen-values larger than $-\gamma$); (ii) an AG method for minimizing a regularized almost-convex function. The first building block is presented in Algorithm~\ref{alg:ncd} and the second building block is proposed in~\citep{DBLP:journals/corr/CarmonDHS16}. For completeness, we also present it in Algorithm~\ref{alg:agac}. The proposed Algorithm (named NEAG) is presented in Algorithm~\ref{alg:neag} and its complexity result is presented in the following theorem.

\begin{theorem}
\label{thm:gnc-a3:iteration}
	With probability at least $1-\delta$, the Algorithm NEAG returns a vector $\xh_k$ such that $\|\nabla f(\xh_k)\|\leq\epsilon$ and $\lambda_{\text{min}}(\nabla^2 f(\xh_k))\geq -\gamma$ with a worse-case iteration complexity of  $\widetilde O\left(\frac{1}{\epsilon\gamma^{3/2}} +\frac{1}{\gamma^{7/2}} +\frac{\gamma^{1/2}}{\epsilon^2}\right)$. 
\end{theorem} 
{\bf Remark:} When $\gamma=\sqrt{\epsilon}$, the iteration complexity of NEAG is $\widetilde O(1/\epsilon^{1.75})$.
\begin{algorithm}[t]
\caption{ \neonp for Negative Curvature Descent: \neonp-NCD$(\x_0, \gamma, c, \delta)$}\label{alg:ncd}
\begin{algorithmic}[1]
\STATE \textbf{Input}:  $\x_0$, $\gamma$, $c>0$, $\delta$
\STATE $\x_1=\x_0$, $\delta' = \delta /(1+12L_2^2\Delta/\gamma^3)$
\FOR{$j=1,2,\ldots,$}
\STATE Compute $\v_j$= \neonp$(f, \x_0, t, \mathcal F, U, \zeta, r)$
\IF{ $\v_j\neq 0$}
\STATE $\x_{j+1} = \x_j - \frac{c\gamma}{L_2}\text{sign}(\v_j^{\top}\nabla f(\x_j))\v_j$
\ELSE
\RETURN $\x_j$
\ENDIF
\ENDFOR
\end{algorithmic}
\end{algorithm}

\begin{algorithm}[t]
\caption{AG for minimizing an almost-convex function: AG-AC$(f,\z_1,\epsilon,\gamma,L_1)$}\label{alg:agac}
\begin{algorithmic}[1]
\FOR{$j=1,2,\ldots$}
\IF{$\|\nabla f(\z_j)\|\leq \epsilon$}
\RETURN $\z_j$
\ENDIF
\STATE Define $g_j(\z)=f(\z)+\gamma\|\z-\z_j\|^2$
\STATE set $\epsilon'=\epsilon\sqrt{\gamma/50(L_1+2\gamma)}$
\STATE $\z_{j+1}=\text{AG-SSC}(g_j,\z_j,\epsilon',L_1,\gamma)$
\ENDFOR
\end{algorithmic}
\end{algorithm}

\subsection{Stochastic Non-Convex Optimization}
As a byproduct, we can also use \neonp~in stochastic non-convex optimization for extracging negative curvature to strengthen first-order stochastic methods for enjoying convergence to a SSP. One can follow the development in~\citep{NEON17} to develop new variants of \neonp-SGD, \neonp-SCSG, \neonp-Natasha with improved convergence, which will be omitted here.

\begin{algorithm}[H]
	\caption{An AG Algorithm for Non-Convex Minimization:  NEAG: $(\x_0, \epsilon, \gamma, \delta)$}\label{alg:neag}
	\begin{algorithmic}[1]
		\STATE \textbf{Input}:  $\x_0$, $\epsilon$, $\gamma$, $\delta$
		\STATE $K:=\lceil1+\Delta\left(\frac{\max(12L_2^2, 2L_1)}{\epsilon_2^3}+\frac{2\sqrt{10}L_2}{\epsilon_1\epsilon_2}\right)\rceil$, $\delta':=\delta/K$, 
		\FOR{$k=1,2,\ldots,$}
		\STATE $\xh_k = \text{\neonp-NCD}(\x_k, \gamma, c, \delta')$
		\IF{$\|\nabla f(\xh_k)\|\leq  \epsilon_1$}
		\RETURN $\xh_k$
		\ELSE
		\STATE  Set $f_k(\x) = f(\x) + L_1\left([\|\x - \xh_j\| - \epsilon_2/L_2]_+\right)^2$
		\STATE $\x_{k+1} = \text{AG-AC}(f_k, \xh_k, \epsilon/2, 3\gamma, 5L_1)$
		\ENDIF
		
		
		\ENDFOR
	\end{algorithmic}
\end{algorithm}

\bibliography{all,ref}

\end{document}